\documentclass[12pt,leqno]{article} 
\usepackage{amsmath,amsfonts,amscd}
\usepackage{pdfsync}
\usepackage{graphicx}
\usepackage{color}
\usepackage[all]{xy}

\usepackage{amssymb,amscd}
\usepackage{graphicx}
\usepackage{color}
\usepackage[all]{xy}
 \usepackage{mathptmx}     

\newcommand{\TO}{\longrightarrow}
\newcommand{\nt}{\noindent}  
\newcommand{\wh}{\widehat}
\newcommand{\per}[1]{\overline{#1}} 
\newcommand{\pro}{\noindent {\it Proof. \/}}
\newcommand{\SF}{\bi{\hbox{\sf S}}{\mathcal F}}
\newcommand{\Hom}{\mathop{\rm Hom \, } \nolimits}
\newcommand{\codim}{\mathop{\rm codim \, } \nolimits}
\newcommand{\Ide}{\mathop{Identity \, } \nolimits}
\newcommand{\Ad}{\mathop{\rm Ad\, } \nolimits}

\newcommand{\supp}{\mathop{\rm supp \, } \nolimits}
\newcommand{\Iso}{\mathop{\rm Iso \, } \nolimits} 
\newcommand{\depth}{\mathop{\rm depth \, } \nolimits}

\newcommand{\rondp}{\hbox{\footnotesize\hbox{$\circ$}}}

\newcommand{\qed}{\hfill $\clubsuit$ \bigskip}
\newcommand{\mf}{M / \mathcal{F}}

\renewcommand\S{\mathbb S}
\newcommand{\refp}[1]{(\ref{#1})}

\newcounter{zatia}
\newcommand{\Zati}{\medskip \setcounter{zatia}{1} (\alph{zatia})   }
\newcommand{\zati}{\medskip \addtocounter{zatia}{1}(\alph{zatia})  }

\newcommand\Om{\Omega}   
\newcommand\om{\omega}

\newcommand\F{\mathcal F}

\newcommand\R{\mathbb R}
\newcommand\HH{\mathbb H}

\newcommand{\bi}[2]{{#1}_{_{#2}}}
\newcommand{\hiru}[3]{{#1}^{^{#2}}{\left( #3 \right)}}
\newcommand{\lau}[4]{{#1}^{^{#2}}_{_{#3}}{\left( #4 \right)}}

\newcommand{\bl}{\begin{lemma}}
\newcommand{\bp}{\begin{proposition}}
\newcommand{\bt}{\begin{theorem}}
\newcommand{\bc}{\begin{corollary}}
\newcommand{\be}{\begin{equation}}
\newcommand{\bee}{\begin{equation*}}

\newtheorem{proposition}{Proposition}[section] 
\newtheorem{lemma}[proposition]{Lemma} 
\newtheorem{corollary}[proposition]{Corollary} 
\newtheorem{theorem}[proposition]{Theorem} 

\newcommand{\el}{\end{lemma}}
\newcommand{\ep}{\end{proposition}}
\newcommand{\et}{\end{theorem}}
\newcommand{\ec}{\end{corollary}}
\newcommand{\ee}{\end{equation}}
\newcommand{\eee}{\end{equation*}}

\newcommand{\bL}{\begin{Lemma}}
\newcommand{\bP}{\begin{Proposition}}
\newcommand{\bT}{\begin{Theorem}}
\newcommand{\bC}{\begin{Corollary}}

\newtheorem{Proposition}{Proposition}[proposition] 
\newtheorem{Lemma}[Proposition]{Lemma} 
  
\newtheorem{Corollary}[Proposition]{Corollary} 
\newtheorem{Theorem}[Proposition]{Theorem}

\newcommand{\eL}{\end{Lemma}}
\newcommand{\eP}{\end{Proposition}}
\newcommand{\eT}{\end{Theorem}}
\newcommand{\eC}{\end{Corollary}}

\newcommand{\pr}[1]{                 
		     \addtocounter{section}{1} 
			\begin{center} \bf \thesection.    #1. \end{center} \setcounter{proposition}{0} }

\newcommand{\prg}[1]{ 
\addtocounter{proposition}{1} {\bf
		\paragraph{\theproposition }  #1. }}

\title{\Large\bf Poincar\'e Duality \\ of the basic intersection cohomology \\of a Killing foliation
\footnote{Mathematical Subjet Classification. 57S15 and 55N33.}}

\author{Martintxo Saralegi-Aranguren\thanks{Laboratoire de Math\'ematiques de Lens, EA~2462 et F\'ed\'eration CNRS Nord- Pas-de-Calais FR~2956. Universit\'e dÕArtois. SP18, rue Jean Souvraz. 62307 Lens Cedex. France.
{\sl saralegi@euler.univ-artois.fr}. }
\\ {\small Universit\'e d'Artois } 
\and   
Robert Wolak\thanks{Instytut Matematyki. Uniwersytet Jagiellonski. 
Stanislawa Lojasiewicza 6, 30\kern 1mm 348
Krakow, Poland.
{\sl robert.wolak@im.uj.edu.pl.}  }
\\ {\small Uniwersytet Jagiellonski}
}

\begin{document}  
\maketitle

\begin{abstract}  We prove that the basic intersection cohomology 
$\lau{\HH}{*}{\per{p}}{\mf}$, where $\F$ is  the singular 
foliation  determined by an isometric action 
of a 
Lie group $G$ on the compact manifold $M$, verifies the  Poincar\'e Duality Property.
 \end{abstract}

\bigskip
 Cohomology theories are one of the basic tools in the study of invariants of topological and geometrical objects. They provide a good environment for the development of obstruction theories. In the case of regular Riemannian foliations basic cohomology theory proved to be of great importance. In particular, for foliations of compact manifolds, it was possible to define a 1-basic cohomology class $\kappa$, the çlvarez class, whose vanishing is responsible for tautness. Moreover, the Poincar\'e duality property holds only in basic twisted cohomology associated to this 1-cohomology class, (see \cite{MR702530} and \cite{MR755165} for the precise statement). 
 
 In the singular case the situation is even more complicated, for isometric actions the top dimensional basic cohomology can vanish and the Poincar\'e duality does not hold, \cite{RoyoPrieto}. Moreover, the standard procedure for the definition of the tautnes class seems not to work. Perhaps one should approach the problem from a different angle, and consider some other cohomology theory.

 We introduced the intersection basic cohomology in \cite{MR2157093} and the examples and results obtained indicate that this cohomology theory is suitable for the study of topology and geometry of singular Riemannian foliations, \cite{MR2403185, MR2551892,MR2262551, MR2968235}. In the present paper we demonstrate that under suitable orientation assumptions the basic intersection cohomology of a Killing foliation satisfies the Poincar\'e duality property.
 
%
%

\medskip

In the sequel $M$  is a connected, second countable, Haussdorff,   
without boundary and smooth 
(of class $C^\infty$) manifold. We also write $G$ for a Lie group.

\newpage

\pr{ Killing Foliations}

 A smooth action $\Phi \colon G \times M 
 \to M$ of the Lie group $G$  is an {\em isometric action} 
when there exists a riemannian metric $\mu$  on $M$ preserved by $G$.
 Moreover, the isometric action $\Phi$ is {\em tame} when the closure of $G$ in 
 $\Iso(M,\mu) $ is  compact.  This is always the case when the manifold $M$ is compact  (cf. \cite{MR1336823}).  
 
 The connected components of the orbits of a tame  action  determine 
  a partition $\mathcal{F}$ 
  on $M$. In fact, this partition is a singular riemannian foliation 
  that we shall call {\em Killing foliation} (cf. \cite{MR932463}).
  
   Notice that $\F$ is also a conical foliation in the sense of \cite{MR2157093}. So, basic intersection cohomology can be used for the study of $\F$. In this work, we prove the Poincar\'e Duality Property of this cohomology.
   
The action $\Phi$ is tame when it is the restriction of   a smooth action $\Phi \colon K \times M 
 \to M$ where $K$ is a compact Lie group containing $G$.
The group $K$ is not unique. We always can choose $K$ in such a 
  way that $G$ is dense in $K$. We shall say that $K$ is a {\em tamer group}. 
  
Since the aim of this work is the study of $\mathcal{F}$ and not the action $\Phi$ itself,  we can consider without loss of generality that the Lie group $G$ is
connected and the action is effective. In this case, $G$ is normal on the tamer group $K$ and the quotient group $K/G$ is commutative
 (see \cite{MR1297165}). This is a key fact for this work.

\pr{Stratification}

 Classifying the points of $M$ following the dimension of the leaves of $\F$
one gets the {\em stratification} $\SF$. It is 
determined by the equivalence relation
$
x \sim y \Leftrightarrow \dim G_{x} = \dim G_{y}.$
The elements  of  $\SF$ are called {\em strata}. The open stratum $R_\F$ is the {\em regular stratum} and the other strata are the {\em singular strata}.

 We fix a basis point $p \in R_\F$ and we put $(G_p)_0 =L$, where  $E_0$  stands for
  the identity component of the Lie group $E$ (the one containing the unity element). This group is generic in the following way.

\bp \label{bat} For each $x\in R_\F$
there exists $k \in K$ with $ \left( G_x \right)_0 = k L k^{-1}$.
Moreover, the choice $x \mapsto k$ can be done locally in a smooth way.
\ep

\pro 
 We consider a  point $x\in R_\F$ and we find an open neighborhood $V \subset R_\F$ of $x$ 
 and a smooth map $f\colon V \to R_\F$ with 
 $
 (G_y)_0 = f(y)(G_x)_0f(y)^{-1},
 $
 for each $y \in V$.
 
 Since the Lie group $K$ is compact we can suppose that  $R_\F = K \bi{\times}{H} \R^m$, where $H$ is a closed subgroup of $O(m)$, and $x= <e,0>$.
We consider a neighborhood $W\subset K$ of $e$. This neighborhood is chosen small enough in order to ensures us the existence of a smooth section $\sigma \colon \gamma (W) \to W$ of the canonical projection $\gamma \colon K\to K/H$. We write 
$\Gamma \colon M =  K\bi{\times}{H} \R^m \to K/H
$
the canonical projection. So, $V =\Gamma^{-1}(W)$ is a neighborhood of $x$ on $R_\F$. Put $f \colon V \to K$ the smooth map defined by $f(y) = \sigma(\Gamma (y))$. A straightforward calculation gives, for each $y=<k,v> = <f(y),v'>$, the equality:
$
G_y = f(y)(G\cap H)_{v'}f(y)^{-1}.
$
On the other hand $G_x = G \cap H$ and $\dim G_x = \dim G_y$ give $\dim (G \cap H) = \dim (G\cap H)_{v'}$ and therefore:
\bee
( G_y )_0= f(y)\left( (G\cap H)_{v'}\right)_0 f(y)^{-1} = 
f(y)\left( G\cap H\right)_0 f(y)^{-1} = f(y)(G_x)_0 f(y)^{-1}.
\eee
This ends the proof.\qed %
 \bigskip
 
 We fix for the rest of the work a  Killing foliation $\F$ given by an effective  tame action $\Phi \colon G \times M \to M$  with $G$ connected.
 We also fix a tamer group $K$. We write $b = \dim G$ and $m =\dim M$. The induced foliation on the regular stratum  $R_\F$ is regular, its dimension will be denoted by $w = \dim \F$.
 
\pr{Presentation of the Poincar\'e Duality Property}

The basic intersection cohomology $\lau{\HH}{*}{\per{p}}{\mf} $, relatively to the perversity $\per{p}$, was introduced in \cite{MR2157093} for the study of conical foliations\footnote{We refer the reader to \cite{MR2968235} for notation and main properties of this notion.}.

 We define the {\em
support} of a perverse form $\om \in \lau{\Pi}{*}{\mathcal{F}}{M}$
as 
$
\supp \om = \overline{\{ x \in M \backslash \Sigma_\mathcal{F} \ / \ \om(x) \not= 0\}},
$
where the closure is taken in $M$.
We denote by
$
\lau{\Om}{*}{\per{q},c}{\mf} =
\left\{ \om \in \lau{\Om}{*}{\per{q}}{\mf} \ \Big/ \ \supp \om
 \hbox{ is compact} \right\}
$
the complex of intersection basic differential forms with compact support relatively to the perversity $\per{q}$. The cohomology $\lau{\HH}{*}{\per{q},c}{\mf}$ of this complex is the {\em
basic intersection cohomology with compact support } of $(M,\mathcal{F})$,
relatively to the perversity $\per{q}$\footnote{We refer the reader to \cite{MR2262551} for notation and main properties of this notion.}.

The goal of this work is to prove that the usual pairing gives the isomorphism 
\bee
\lau{\HH}{*}{\per{p}}{\mf} \cong \lau{\HH}{m-w -*}{\per{q},c}{\mf},
\eee
where $\per{p}$ and $\per{q}$ are complementary perversities, that is, $\per{p} + \per{q} = \per{t}$, with $\per{t}(S) = \codim_M \F  \!  -  \! \codim_S \F_S \! -   \!  2$ where $S$ is a singular stratum and $\F_S $ the restriction of $\F$ to $S$.

\pr{Twisted product}
The elementary pieces on $M$ are the twisted products. We find in \cite[Proposition 5]{MR2968235} the computation of   their basic intersection cohomology. We present here the compact support version of this result.

\bp \label{supco}
$
\lau{\HH}{*}{\per{q},c}{K \bi{\times}{H} N / \mathcal{W} } = \left( \hiru{H}{*} {K / \mathcal{K} }  \otimes \lau{\HH}{*}{\per{q},c}{N /  \mathcal{N}  }\right)^{H/H_{0}}.
$
\ep
\pro It suffices to follow  \cite[Proposition 5]{MR2968235} taking into the account that, given a differential form $\om$ on 
$K \bi{\times}{H} R_{\mathcal{W}}$, we have:
\bee
\Pi^* \om \in \lau{\Om}{*}{\per{q},c}{K \times N /\mathcal{E}\times
 \mathcal{N}} \Longleftrightarrow \om \in \lau{\Om}{*}{\per{q},c}{K \bi{\times}{H}N / \mathcal{W}}.
\eee
This comes from the fact that $\Pi$ is an onto map and that the Lie groups $K$ and $H$ are compact. \qed

\pr{Tangent volume form}

In order to construct the pairing giving the Poincar\'e Duality we need to introduce a particular tangent volume form of the orbits of $\Phi$.

We  fix a $K$-invariant metric $\nu$ on $\mathfrak{k}$, the Lie algebra of $K$ which exists since $K$ is compact. Consider $ \{ u_1, \ldots , u_f\}$ an orthonormal basis of $\mathfrak{k}$   where $\{u_1, \cdots , u_b\}$ is a basis of  
 $\mathfrak{g} $ and  $\{u_1, \cdots , u_w\}$ is a basis of  
 $\mathfrak{l}^{\bot} $. Here,  $\mathfrak{g}$ (resp. $\mathfrak{l}$) denotes the Lie algebra of  $G$ (resp. $L$).
We put  $\tau = u_1^* \wedge \cdots \wedge u_w^*$ 
 the associated volume form of $\mathfrak{l}^{\bot}$. 

  We write $V_u$ the fundamental vector field  on $M$ associated to $u \in \mathfrak{g}$. A {\em tangent volume form} of $\Phi$ is a $G$-invariant differential form
$\eta \in \lau{\Pi}{w}{\mathcal{F}}{M}$ verifying:
\be\label{def}
\eta(V_{ u_1}(x),\ldots,V_{ u_{w}}(x)) =
\tau \left( \Ad(\ell^{-1}) \cdot v_1,\ldots, \Ad(\ell^{-1}) \cdot v_{w}\right),
\ee
where  $\{v_1, \ldots , v_{w}\} \subset \mathfrak{g}$, $x \in R_\F$ and $G_x = \ell L \ell^{-1}$.

We prove the existence of a tangent volume form under   a suitable orientation conditions on the manifold and on the foliation. The tame action $\Phi $ is {\em orientable} if
 \begin{itemize}
 \item[(i)] the manifold $M$ is orientable, and
 \item[(ii)] the adjoint action $\Ad \colon N_K(L) \times \mathfrak{l} \to \mathfrak{l}$  is orientation preserving.
 \end{itemize}

\prg {\bf Remarks}

\Zati Condition (ii) does not depend on the choice of $L$.

\zati Since the group $K$ preserves the orientation of  $\mathfrak{g}$ then condition (ii) is equivalent to
 \begin{itemize}
 \item[]
  \begin{itemize}
  \item[]
   \begin{itemize}
 \item[(ii')] the adjoint action $\Ad \colon N_K(L) \times \mathfrak{l}^\bot \to \mathfrak{l}^\bot$  is orientation preserving.
 \end{itemize}
 \end{itemize}
  \end{itemize}
  
\zati Condition (ii) is verified when $G$ is abelian or when $\frak{l}=0$, that is, when  $\dim \F = \dim G$.


\bp
\label{etatang} If the action $\Phi$ is orientable, then there exists a $K$-invariant tangent volume form
$\eta$
of $\Phi$
verifying:

\Zati  
%
For each $\om \in \lau{\Om}{m -w - 1}{\per{t}}{\mf}$ the
product $\om \wedge d\eta$ is 0.

\zati For each $\om \in \lau{\Om}{m -w}{\per{t},c}{\mf}$ the
integral $\displaystyle{\int_{R_{\mathcal{F}}} \omega \wedge \eta }$ is finite.

\zati For each $\om \in \lau{\Om}{m -w-1}{\per{t},c}{\mf}$ the
integral $\displaystyle{\int_{R_{\mathcal{F}}} d(\omega \wedge \eta) }$
is 0.
\ep
\pro Firstly, we  prove the following statement by induction
on $\depth \SF$:
\begin{quote}
    ``There exists a $K$-invariant differential form
    $\overline\eta \in \lau{\Pi}{w}{\F \times \mathcal{I}}{M\times
    [0,1[^p}$ verifying 
    \be\label{defBis}
\overline\eta((V_{ v_1}(x),0),\ldots, (V_{ v_{w}}(x),0) ) =\tau( \Ad(\ell^{-1}) \cdot v_1,\ldots,\Ad(\ell^{-1}) \cdot  v_{w}),
\ee
    where  $\{v_1, \ldots , v_{w}\} \subset \mathfrak{g}$, $x \in R_\F$ and $(G_x)_0 = \ell L\ell^{-1}$ with $\ell\in 
    K$.''
    \end{quote}
Here, $\mathcal{I}$ denotes the pointwise foliation of $[ 0,1[^p$.    The existence of $\eta$ is proven by taking $p=0$.

\begin{center} {\em First case: $\depth \SF = 0$.}\end{center}

We have $\lau{\Pi}{*}{\F}{M\times
    [0,1[^p} = \lau{\Om}{*}{}{(M\times
    [0,1[^p)/(\F\times \mathcal{I})}$. Since the foliation $\F$ is  $K$-invariant then it suffices to define $\overline\eta$ on $T(\F \times \mathcal{I})$. In fact, this restriction is given by \refp{defBis}. It remains to prove that $\overline\eta$ is well-defined, smooth on $T(\F \times \mathcal{I})$ and $K$-invariant. Let us see that.

\begin{itemize}
\item {\em The definition \refp{defBis} does not depend on $\ell$}. Let us consider $\ell' \in K$ with $(G_x)_0 =\ell' L \ell'^{-1}$. 
Then $  \ell'^{-1} \ell  \in N_K(L)$. This gives
\bee
\tau(\Ad(\ell'^{-1}) \cdot v_1,\ldots,\Ad(\ell'^{-1}) \cdot v_w) = \tau(\Ad ( \ell'^{-1} \ell )\Ad(\ell^{-1}) \cdot v_1),\ldots, \Ad ( \ell'^{-1} \ell  )\Ad(\ell^{-1})   ) \cdot v_w).
\eee
Since  the element $\Ad ( \ell' \ell^{-1} ) $ preserves the metric $\nu$ and the orientation of $\mathfrak{l}^\bot$ (see (ii')) then we get
\bee
\tau(\Ad(\ell'^{-1}) \cdot v_1,\ldots,\Ad(\ell'^{-1}) \cdot v_w) = \tau(\Ad(\ell^{-1}) \cdot  v_1,\ldots,   \Ad(\ell^{-1}) \cdot v_w).
\eee

\item   {\em The definition \refp{defBis} is smooth}.
Consider $x \in M$. From Proposition \ref{bat} we know that there exists  a neighborhood $V\subset M$ and a smooth map $f \colon V \to K$ such that $(G_y)_0 = f(y) Lf(y)^{-1}$ for each $y \in V$. In this neighborhood we have
$
\overline\eta((V_{v_1}(y),0)\ldots, (V_{v_w}(y),0)) =\tau(\Ad(f(y)^{-1})\cdot v_1),\ldots,\Ad(f(y))^{-1})\cdot v_w),
$
which is smooth.

\item {\em The form $\overline\eta$ is $K$-invariant}. If $k \in K$ we get $(G_{k \cdot x})_0 = k \ell L \ell^{-1}k^{-1}$ and 
\begin{eqnarray*}
(k^*\overline\eta)((V_{v_1}(x),0),\ldots,(V_{v_w}(x),0)) =
\overline\eta(k_*(V_{v_1}(x),0),\ldots,k_*(V_{v_w}(x),0))&=& \\
\overline\eta((V_{\Ad (k) \cdot v_1}(k \cdot x),0)\ldots, (V_{\Ad (k) \cdot v_w}(k \cdot x),0)) 
= \tau(\Ad((k\ell)^{-1}) \Ad (k) \cdot v_1,\ldots,\Ad((k\ell)^{-1}) \Ad (k) \cdot v_w) &=& \\
 \tau(\Ad(\ell^{-1}) \cdot v_1,\ldots,\Ad(\ell^{-1})\cdot v_w) =\overline \eta((V_{v_1}(x),0)\ldots,(V_{v_w}(x),0)).
\end{eqnarray*}
\end{itemize}

\begin{center} {\em Second case: $\depth \SF > 0$.}\end{center}

By induction hypothesis there exists  a $K$-invariant differential
form $\overline\eta_{0}
\in \lau{\Pi}{w}{\wh\F}{\wh{M}\times
    [0,1[^p}$ verifying \refp{defBis}. Associated to
the Molino's blow up (cf. \cite[5.2]{MR2968235}) we have the $
K$-equivariant imbedding
$
 \sigma \colon M\backslash
S_{_{min}} \to \mathcal{L}^{-1}(M\backslash
S_{_{min}}),
$
defined by
$
\sigma(z) = (z,1).
    $
    The differential form $\overline\eta = (\sigma\times \hbox{identity}_{[0,1[^p})^{*}\eta_{0}$
    belongs to the complex $\hiru{\Om}{w}{R_{\F}\times
    [0,1[^{p}}$, is $K$-invariant and verifies \refp{defBis}. It remains to prove that
    $\overline\eta \in \lau{\Pi}{\ell}{\F \times \mathcal{I}}{M \times
    [0,1[^p}$, which is a local property.
    
So, we  can consider that $M$ is a tubular
 neighborhood  $T$ of a singular stratum of $\SF$ and  prove
 $(\nabla\times \hbox{identity}_{[0,1[^{p+1}})^{*} \overline\eta \in
\lau{\Pi}{w}{\F \times \mathcal{I}}{D \times [0,1[^{p+1}} $ (cf. \cite[3.1.1 (e)]{MR2262551}).
This is the case since  $ \sigma \rondp \nabla \colon D \times ]0,1[ \to D \times ]-1,1[$ is just the
 inclusion and $\overline\eta_{0} \in \lau{\Pi}{w}{\F \times \mathcal{I}}{D\times ]-1,1[ \times
    [0,1[^{p}}$.
    
\medskip

We prove now the (a)-(c) items.

\Zati 
%
It suffices to prove this property on the regular stratum $R_\F$. In other words, we can suppose that $\depth \SF =0$.
 Since the question is a local one, it is enough to prove $\om \wedge d\eta= 0$ on the open subset $V \subset M$ (cf. proof of Proposition \ref{bat}).
 
 For each $y \in V$ we have $(G_y)_0 =f(y)L f(y)^{-1} $. Then, 
 $\{ V_{f(y)\cdot u_1}(x), \ldots, V_{f(y)\cdot u_w}(x)\}$ is a basis of $T_y G(y)$. 
  For degree reasons it suffices to prove that we have $i_{V_{f(y)\cdot u_1}(x)} \cdots i_{V_{f(y)\cdot u_w}(x)}
(\om \wedge d\eta)= 0$. Since $\om$ is a basic form and $\eta$ is a $K$-invariant form, we can write
\begin{eqnarray*}
i_{V_{f(y)\cdot u_1}(y)} \cdots i_{V_{f(y)\cdot u_w}(y)}
(\om \wedge d\eta)
&=&
 (-1)^{w }\om \wedge i_{V_{f(y)\cdot u_1}(y)} \cdots i_{V_{f(y)\cdot u_w}(y)}d\eta =\om \wedge d \left( i_{V_{f(y)\cdot u_1}(y)} \cdots i_{V_{f(y)\cdot u_w}(y)} \eta\right)  \\
&= &\om \wedge d \left(\eta \left(V_{f(y)\cdot u_1}(y) , \ldots  , V_{f(y)\cdot u_w}(y)\right)\right)=
\om \wedge d \left(\tau \left(u_1 ,\ldots ,  u_w\right)  \right)  \\ &=& \om \wedge d  1 = 0. \end{eqnarray*}

\zati  and \zati. Notice that the integrals make sense since $M$ is an oriented manifold and $R_\F$ is an open subset of it. Now, the proof is the same as in in \cite[Lemma 4.3.2]{MR2262551}.
\qed

\pr{The pairing} In Section 9 we  prove the Poincar\'e Duality
Property:
$
\lau{\HH}{*}{\per{p}}{\mf} \cong \lau{\HH}{m-w -*}{\per{q},c}{\mf},
$ when  $\Phi$ is orientable and the   two  perversities $\per{p}$ and
$\per{q}$ are complementary.
This isomorphism comes from the pairing $P_M$ constructed from the above  tangent volume form $\eta$ in the following way:
$$
\xymatrix{
P_M\colon 
\lau{\Om}{*}{\per{p}}{\mf} \times
\lau{\Om}{m-w-*}{\per{q},c}{\mf} \longrightarrow \R & \therefore &(\alpha,\beta) \leadsto   \displaystyle \int_{R_\mathcal{F}}
\alpha \wedge \beta \wedge \eta.
}
$$
Proposition
\ref{etatang} implies that this operator is well defined and that it induces the {\em pairing}
\bee
P_M
 \colon
\lau{\HH}{*}{\per{p}}{\mf} \times
\lau{\HH}{m-w-*}{\per{q},c}{\mf} \TO \R,
\eee
defined by
$
\displaystyle P_M ([\alpha],[\beta]) = P_M (\alpha,\beta).
$
The Poincar\'e Duality Property asserts  that $P_M$ is a
non degenerate pairing, that is, the operator
\bee
P_M
 \colon
\lau{\HH}{*}{\per{p}}{\mf} \TO \Hom
\left(\lau{\HH}{m-w-*}{\per{q},c}{\mf} , \R
\right)
\eee
defined by
$
\displaystyle 
P_M ([\alpha])([\beta]) = \int_{R_\mathcal{F}} \alpha \wedge \beta \wedge \eta
$
is an isomorphism.

\pr{Twisted product and Poincar\'e Duality}

We first get the Poincar\'e Duality Property in the framework of twisted products.

\bp
\label{modelDP}
Consider a twisted product $ K \times_{H } N$, where $N$ is an orientable manifold and the action $\Theta \colon H \times N \to N$ is effective. Let us suppose that the associated action $\Phi \colon G \times  (K \times_{H } N )\to (K \times_{H } N)$ is orientable. Then:

\Zati The action $ \Theta \colon (G \cap H)_0 \times N \to N$ is orientable.

\zati Si $
(N,\mathcal{N})$ verifies the Poincar\'e Duality Property  then $
(K \times_{H }N, \mathcal{W})$ verifies the
Poincar\'e Duality Property.
\ep

\pro  Recall that the action $ \Theta \colon (G \cap H)_0 \times N \to N$ is a tame action relatively to the closure $H'$ of $(G \cap H)_0$ on $H$. 

For each $<k,v> \in K \times_{H }N$ we have $G_{<k,v>} = k  \ (G\cap H)_v \ k^{-1}$.
Then $((G\cap H)_v)_0 =\ell L \ell^{-1}$ and
$(G_{<k,v>})_0 = k\ell L \ell^{-1}k^{-1}$ for some $\ell \in H'$.
 In particular, we can take the same $L$ for both actions $\Phi$ and $\Theta$. 
Property (a) comes now from the inclusion $N_{H'} (L) \subset N_K (L)$.
We prove (b) in several steps.
 
 \medskip

1. {\em Tangent volume forms of $\Phi$ and $\Theta $.} We follow the notations of \cite[4.1]{MR2968235}. In particular, we consider
\begin{equation*}
B = \left\{ u_1, \ldots u_a,u_{a+1}, \ldots , u_w, u_{w+1}, \ldots  , u_{b}, u_{b+1},\ldots ,
u_c,u_{c+1},\ldots,u_f\right\}
\end{equation*}
an orthonormal basis of $\mathfrak{k}$ with
 $\left\{ u_1, \ldots u_b\right\}$ basis of  $\mathfrak{g}$,   
$\left\{ u_{a+1}, \ldots u_c\right\}$ basis of the Lie algebra $\mathfrak{h}$ of $H$ and $\left\{ u_{w+1}, \ldots u_b\right\}$ basis of $\mathfrak{l}$.

Consider a tangent volume form $\eta$ of $\Phi \colon G \times M \to M$ (resp. $\eta_0$ of $\Theta \colon (G \cap H)_0 \times N \to N $) 
associated to the metric $\nu$ (resp. of $\nu_{|{\mathfrak{g} \cap \mathfrak{h}}}$).  Recall that $\tau = u_1^* \wedge \cdots \wedge u_w^*$ and 
$\tau_0 = u_{a+1}^* \wedge \cdots \wedge u_w^*$.
We prove that
\be\label{eta}
 (-1)^{a(b-a)}\gamma_{a+1} \wedge \cdots \wedge \gamma_b\wedge \Pi^* \eta  =\gamma_1 \wedge \cdots \wedge \gamma_b\wedge \eta_0 \ \ \hbox{ on }  \mathcal{K} \times \mathcal{N}.
\ee

The leaf of $\mathcal{K} \times \mathcal{N}$ at the point $(k,v) \in K \times N$ is generated by
\bee
 \mathfrak{B} =\left\{(X_{u_1}(k),0), \ldots , (X_{ u_b}(k),0), (0,W_{\Ad(\ell) (u_{a+1})}(v)), \ldots , (0,W_{\Ad(\ell)(u_w)}(v))\right\}
\eee  
 (cf. \cite[Proposition 5 $\langle  v \rangle $]{MR2968235}).The RHT of \refp{eta} applied to $\mathfrak{B} $ gives, 
\begin{eqnarray*}
\eta_0(W_{\Ad(\ell) (u_{a+1})}(v)), \ldots , W_{\Ad(\ell)(u_w)}(v) =
\tau_0(u_{a+1}, \ldots , u_w) =
(u_{a+1}^* \wedge \cdots \wedge u_w^*) (u_{a+1}, \ldots , u_w) =1.
\end{eqnarray*}
The LHT of \refp{eta} applied to $\mathfrak{B} $ gives, using the fact that $\Pi_*(X_u(k),-W_u(v))=0$ if $u \in \mathfrak{g} \cap \mathfrak{h}$:
\begin{eqnarray*}
\Pi^*\eta \left((X_{u_1}(k),0), \ldots , X_{ u_a}(k),0), (0,W_{ \Ad(\ell)(u_{a+1})}(v)), \ldots , (0,W_{\Ad(\ell)(u_w)}(v))\right) &= &\\
 \eta \left( \Pi_*(X_{u_1}(k),0), \ldots , \Pi_*(X_{ u_a}(k),0), \Pi_*(X_{ \Ad(\ell)(u_{a+1})}(k),0), \ldots , \Pi_*(X_{\Ad(\ell)(u_w)}(k),0)\right) &=&\\
 \eta \left( \Pi_*(X^{\Ad(k)(u_1)}(k),0), \ldots , \Pi_*(X^{\Ad(k)( u_a)}(k),0), \Pi_*(X^{\Ad(k\ell)( u_{a+1})}(k),0), \ldots , \Pi_*(X^{\Ad(k\ell)(u_w)}(k),0)\right)&=& \\
\eta \left( V_{\Ad(k)(u_1)}(<k,v>) \ldots , V_{\Ad(k)(u_a)}(<k,v>),  V_{ \Ad(k\ell)(u_{a+1})}(<k,v>), \ldots , V_{\Ad(k\ell)(u_w)}(<k,v>)\right)&=&\\
\tau \left(\Ad(\ell^{-1})(u_1),\ldots , \Ad(\ell^{-1})(u_a), u_{a+1}, \ldots , u_w\right)=
 (u_1^* \wedge \cdots \wedge u_a^*) \left(\Ad(\ell^{-1})(u_1)\ldots , \Ad(\ell^{-1})(u_a)\right).
 \end{eqnarray*}
 Here, $V_u$ denotes the fundamental vector field
 of the action $\Phi$ associated to $u\in \mathfrak{g}$ . Since $\ell \in H'$,  the closure of $(G \cap H)_0$ on $H$, then $\Ad(\ell)$ preserves $\mathfrak{g}$ and $\mathfrak{h}$. Connectedness of $H'$ gives that the operator $\Ad(\ell) \colon (\mathfrak{g} \cap \mathfrak{h} )^{\bot_\mathfrak{g}} \to ( \mathfrak{g} \cap \mathfrak{h})^{\bot_\mathfrak{g}}$ is an orthogonal map preserving orientation. So
 \begin{eqnarray*}
  (u_1^* \wedge \cdots \wedge u_a^*) \left(\Ad(\ell^{-1})(u_1)\ldots , \Ad(\ell^{-1})(u_a)\right) =
\det (\Ad(\ell^{-1}) )  (u_1^* \wedge \cdots \wedge u_a^*) \left(u_1,\ldots , u_a\right)&=&
1,
\end{eqnarray*}
We obtain  \refp{eta}.
 
2. {\em Some maps}.
Consider now two complementary perversities $\per{p}$ and $\per{q}$ on $N$. We have $\dim N =m+c-a-f$ and $\dim \mathcal{N}=w-a$, where $w = \dim \mathcal{W}$ and $m = \dim K \bi{\times}{H} {N}$. 
By
hypothesis, the pairing
$
P_{N} \colon 
\lau{\HH}{*}{\per{p}}{N/\mathcal{N}}\times
\lau{\HH}{m+c-w-f - *}{\per{q},c}{N/\mathcal{N}}\to \R
$
is  non degenerate.  On the other hand, it is clear that the pairing
$
P \colon
\hiru{H}{*}{K/\mathcal{E}}
\times
\hiru{H}{f-c-*}{K/\mathcal{E}}
\TO \R,
$
defined by
$
P([\xi],
[\chi]) =
{\displaystyle \int_{K}  \xi
\wedge \chi \wedge \gamma_1 \wedge \cdots \wedge \gamma_c}
$
is non degenerate (cf. \cite[4.1]{MR2968235}. So, the first row of the below diagram 3. is non degenerate.

Write also $\per{p}$ and $\per{q}
$ the associated perversities on $K \times_H N$, which also are two complementary perversities.
Recall that the isomorphisms 
\begin{eqnarray*}
&& \nabla  \colon \left( \hiru{H}{*}{K/\mathcal{E}} \otimes \lau{\HH}{*}{\per{p}}{N/\mathcal{N}}\right)^{H/H_0} \TO \lau{\HH}{*}{\per{p}}{K \times_H N}, \\
&& \nabla  \colon \left( \hiru{H}{*}{K/\mathcal{E}} \otimes \lau{\HH}{*}{\per{q},c}{N/\mathcal{N}}\right)^{H/H_0} \TO \lau{\HH}{*}{\per{q},c}{K \times_H N} 
\end{eqnarray*} 
 are characterized by
 \be\label{pinabla}
\displaystyle  \Pi^* \nabla ([\xi] \otimes  [\alpha]) = \left[\xi \wedge \left(\alpha + \sum_{b < i_1 < \cdots < i_l\leq c} \ (-1)^\ell \gamma_{i_1} 
\wedge \dots 
\wedge \gamma_{i_l} \wedge 
(i_{W{_{i_l}} }\cdots i_{W{_{i_1}}}\alpha)\right)\right]
 \ee
(cf. \cite[Proposition 5]{MR2968235} and Proposition \ref{supco}).

\medskip

3. {\em A diagram}.
Let us consider the following diagram
\bee
\begin{CD}
\left(\hiru{H}{*}{K/\mathcal{E}}
\otimes 
\lau{\HH}{*}{\per{p}}{N/\mathcal{W}}\right)^{H/H_0}
\! \times \!
\left(\hiru{H}{f-c-*}{K/\mathcal{E}}
\otimes 
\lau{\HH}{m+c-w-f-*}{\per{q},c}{N/\mathcal{N}}\right)^{H/H_0}
@>P \otimes P_{N}>> \R\\
    @V \nabla\times \nabla  VV @V\Ide VV\\
\lau{\HH}{*}{\per{p}}{K \times_H N}\times
\lau{\HH}{m-w -*}{\per{q},c}{K \times_H N}
@>P_{K \times_H N} >> \R.\\
    \end{CD}
    \eee
We end the proof if we show that this diagram  commutates up to a constant. Write $\nabla ([\xi] \otimes  [\alpha])  =[\alpha \bullet \beta]$. We have,
\bee
P_{K\times_H N} (\nabla\times \nabla) ([\xi]
\otimes [\alpha],[\chi]
\otimes [\beta]) 
=
\displaystyle \int_{K  \bi{\times}{H} R_\mathcal{N}}
\xi \bullet \alpha
\wedge
\chi\bullet \beta
 \wedge \eta.
 \eee
Recall that 
we have  denoted by
$\{W_{a+1},\ldots , W_c\}$ the
fundamental vector fields of the action $\Theta \colon H \times N \to 
N $ associated to the basis $\{u_{a+1}, \ldots , u_{c}\}$. We write $\{\zeta_{a+1}, \ldots , \zeta_c\}$ the associated dual forms relatively to an $H$-invariant riemannian metric on $R_{\mathcal{N}}/\mathcal{N}$.  So, $\frac{1}{2^{c-a}}(\gamma_{a+1} + \zeta_{a+1}) \wedge \cdots \wedge(\gamma_{c} + \zeta_{c})$ is a differential form of $K \times N$ giving a volume form on each fiber of $\Pi$. Thus
\begin{eqnarray*}
P_{K\times_H N} (\nabla\times \nabla) ([\xi]
\otimes [\alpha],[\chi]
\otimes [\beta]) 
&=&
\frac{1}{2^{c-a}}\displaystyle \int_{K\times R_\mathcal{N}}
\Pi^* (\xi \bullet \alpha
\wedge
\chi\bullet \beta
 \wedge \eta ) \wedge (\gamma_{a+1} + \zeta_{a+1}) \wedge \cdots \wedge(\gamma_{c} +\zeta_{c}) \\
 &\stackrel{\refp{pinabla}}{=}&
 \frac{1}{2^{c-a}}\displaystyle \int_{K\times R_\mathcal{N}}
\xi \wedge \alpha
\wedge
\chi \wedge \beta
 \wedge \Pi^* \eta \wedge (\gamma_{a+1} + \zeta_{a+1}) \wedge \cdots \wedge(\gamma_{c} +\zeta_{c}) 
 \end{eqnarray*}
  for degree reasons. By reordering
\begin{eqnarray*}
P_{K\times_H N} (\nabla\times \nabla) ([\xi]
\otimes [\alpha],[\chi]
\otimes [\beta]) 
=
\frac{\pm 1}{2^{c-a}}\int_{K  \times R_\mathcal{N}} \alpha \wedge \beta \wedge 
\gamma_{a+1}\wedge \cdots \wedge \gamma_{f} \wedge \Pi^* \eta &\stackrel{\refp{eta}}{=}&\\[,3cm]
\frac{\pm 1}{2^{c-a}}\int_{K  \times R_\mathcal{N}} \alpha \wedge \beta \wedge 
\gamma_{1}\wedge \cdots \wedge \gamma_{f} \wedge  \eta_0= \frac{\pm 1}{2^{c-a}}\int_{K } \gamma_{1}\wedge \cdots \wedge \gamma_{f}  \cdot  
\int_{ R_\mathcal{N}}\alpha \wedge \beta \wedge 
 \eta_0
&=&\\[,3cm]
 \frac{\pm 1}{2^{c-a}}P([\xi],[\chi])  \cdot
P_{N}([\alpha],[\beta]) =
\frac{\pm 1}{2^{c-a}}(P \otimes P_{N})([\xi] \otimes [\alpha],[\chi] \otimes [\beta]).  &&
\end{eqnarray*}
we obtain the commutativity. \qed

\pr{Tubular neighborhoods and Poincar\'e Duality}

Consider
$(T,\tau,S,\R^n)$  a $K$-invariant tubular neighborhood of a
singular stratum $S$.  The restriction of the action $\Phi \colon G \times M \to M$ to $T$ is orientable. 
Put $(\R^n,\F_{\R^n})$ the slice of the tubular
neighborhood. That is, $\R^n$ is identified with a fiber $\tau^{-1}(x), x \in S$ and the foliation $\F_{\R^n}$ is defined by a tame action $\Theta \colon G_x\times \R^n \to \R^n$. It is also an orientable action since $(G_y)_x= G_y$ for each $y \in \tau^{-1}(x)$.

\bp \label{tubDP} 
 If  the
slice verifies the Poincar\'e Duality
Property then the tube  also verifies
the Poincar\'e Duality Property.
 \ep 
 \pro
 The proof is the same of that of \cite[Propositon 6]{MR2968235} by
considering the following statement
\bee
\mathfrak{A}(M,\mathcal{F}) = 
\left\{
\begin{array}{c}
\hbox{
``The pairing
    $
    P_{T} \colon \lau{\HH}{*}{\per{p}}{T/\mathcal{F}} \TO
\lau{\HH}{m-w - *}{\per{q}}{T/\mathcal{F}} 
    $
    is non degenerate,} \\ \hbox{for any two complementary perversities $\per{p}$ and $\per{q}$."}
\end{array}
\right.
\eee
and replacing 4.3 (b) by  Proposition \ref{modelDP}.\qed

\pr{Poincar\'e Duality}

The main result of this work is the following

 \bt
\label{TPC} 
The basic intersection cohomology associated to there Killing foliation $\mathcal{F}$ determined  by an orientable 
action of a Lie group verifies the Poincar\'e Duality Property. 
\et
 \pro 
The proof of the Theorem is the same of that of  
\cite[Theorem 1]{MR2968235} by considering the statement 
\bee
\mathfrak{A}(M,\mathcal{F}) = 
\left\{
\begin{array}{c}
\hbox{``The pairing 
$P_M \colon\lau{\HH}{*}{\per{p}}{\mf} \times \lau{\HH}{m-w-*}{\per{q},c}{\mf} \TO \R$, is non degenerate }\\
\hbox{
for any two complementary perversities $\per{p}$ and $\per{q}$,"}
\end{array}
\right.
\eee
\nt and the following remarks:

\begin{itemize}
\item The result for the regular foliation comes form
 $\cite{MR865667}$. 

\item The Mayer-Vietors sequence used are those of \cite[3.6, 3.7]{MR2262551}.

\item The reference \cite[Proposition 6]{MR2968235} is replaced by Proposition
\ref{tubDP}.
\end{itemize}
So,
 we reduce the problem to prove
$
\mathfrak{A}(\R^{m_2},\mathcal{F}_{\R^{m_2}})$, that is, to prove that the
pairing
$$
P_{\R^{m_2}} \colon \lau{\HH}{*}{\per{p}}{\R^{m_2}/\mathcal{F}_{\R^{m_2}}}  \times
\lau{\HH}{m_2-w_2 -*}{\per{q},c}{\R^{m_2}/\mathcal{F}_{\R^{m_2}}}  \TO \R,
$$
is non degenerate. Here, $\mathcal{F}_{\R^{m_2}}$ is the foliation given by an orientable orthogonal (and therefore tame) action  $\Lambda_x \colon G_x \times \R^{m_2} \to \R^{m_2}$ having the origin as the unique fixed point.  We put $w_2 = \dim \mathcal{F}_{\R^{m_2}}$. We write $\mathcal{G}$ induced  foliation on $\S^{m_2-1}$. By the induction on the depth of $\SF$ , we have $\mathfrak{A}(\S^{m_2-1},\mathcal{G}) $.

From \cite[Proposition 4]{MR2968235}  we have
\begin{eqnarray}
\label{cono1}
\lau{\HH}{i}{\per{p}}{\R^{m_2}/\mathcal{F}_{\R^{m_2}}} &=&
\left\{
\begin{array}{cl}
\lau{\HH}{i}{\per{p}}{\S^{m_2-1}/\mathcal{G}}  & \hbox{if } i \leq
\per{p}(\vartheta)\\[,3cm]
0 & \hbox{if } i \geq
\per{p}(\vartheta)+1.
\end{array}
\right.
\end{eqnarray}
From \cite[Proposition  3.7.2]{MR2262551} and the fact that $\per{p}$ and $\per{q}$ are complementary perversities on $\R^{m_2}$; id set, 
$
\per{p}(\{ \vartheta\}) + \per{q}(\{ \vartheta\})
=
\per{t}(\{ \vartheta\})
=
m_2  - w_2  - 2,
$
we have
\begin{eqnarray}
\label{cono2}
\lau{\HH}{m_2 -w_2 -i}{\per{q},c}{\R^{m_2}/\mathcal{F}_{\R^{m_2}}} &=&
\left\{
\begin{array}{cl}
\lau{\HH}{m_2 - w_2 -i-1}{\per{q}}{\S^{m_2-1}/\mathcal{G}}  &
\hbox{if } i \leq
\per{p}(\{ \vartheta\})  \\[,3cm]
0 & \hbox{if } i\geq
\per{p}(\{ \vartheta\})   + 1 .
\end{array}
\right.
\end{eqnarray}
Now, $\mathfrak{A}(\R^{m_2},\mathcal{F}_{\R^{m_2}})$  comes from $\mathfrak{A}(\S^{m_2-1},\mathcal{G}) $ and these two facts
\begin{itemize}
\item[(i)]  The pairing $P_{\R^{m_2}}$ becomes the
pairing $P_{\S^{m_2-1}}$ through the
isomorphism induced by \refp{cono1} and \refp{cono2}.
\begin{quote}
   Notice first that a tangent volume form $\eta$ of  $\mathcal{G}$ is also a tangent volume form  for
$\mathcal{F}_{\R^{m_2}}$.  The operator
     $
     \aleph \colon
     \lau{\HH}{*}{\per{p}}{\S^{m_2-1}/\mathcal{G}_S}
     \to
     \lau{\HH}{*}{\per{p}}{c\S^{m_2-1}/
     c\mathcal{G}_S}
     $
     defining \refp{cono1} is
     $
     \aleph([\alpha]) = [\alpha];
     $  the operator
     $
     \aleph' \colon
     \lau{\HH}{*}{\per{q}}{\S^{m_2-1}/\mathcal{G}_S}
     \to
     \lau{\HH}{*}{\per{q},c}{c\S^{m_2-1}/
     \mathcal{G}_S}
     $
     defining \refp{cono2} is $
     \aleph'([\beta]) = [g\,
     dt \wedge \beta].
     $
Now, for
$[\alpha] \in \lau{\HH}{i}{\per{p}}{\S^{m_2-1}/\mathcal{G}_S}$
and
$
[\beta] \in
\lau{\HH}{m_2 - 1 - w_2-i}{\per{q},c}{\S^{m_2-1}/\mathcal{G}_S}
$
we have
\bee
P_{\R^{m_2}} (\aleph_{23}[\alpha], \aleph'[\beta])=
\int_{ R_{\mathcal{G}} \times ]0,1[} \alpha
\wedge  g \,
\wedge dt \wedge \beta \wedge \eta =
\left(\int_{ R_{\mathcal{G}}} \alpha \wedge \beta \wedge \eta\right)
\left(\int_{0}^1 g dt\right)
=
P_{\S^{m_2-1}} ([\alpha], [\beta]).
\eee

\end{quote}
\item[(iii)] The perversities $\per{p}$ and $\per{q}$ are complementary on
$\S^{m_2-1}$.
\begin{quote}
We have, for any stratum $S \in \hbox{\sf S}_{{\mathcal G}}$
the equalities
$\per{p}(S) + \per{q}(S) =\per{p}(S \times ]0,1[) + \per{q}( S \times ]0,1[)=
\per{t}(S \times ]0,1[)
= \codim_{\R^{m_2}} \mathcal{F}_{\R^{m_2}} -
\codim_{S \times ]0,1[} (
\mathcal{G}_S \times \mathcal{I}) -2 = \codim_{\S^{m_2-1}}
\mathcal{G} -\codim_S \mathcal{G}_S -2 = \per{t}(S)$.
\end{quote}
\end{itemize}
Hau amaiera da. 
\qed

\bibliographystyle{abbrv}

\end{document}